\magnification=1200
\font\titlefont=cmcsc10 at 12pt
\hyphenation{moduli}

%
%
\catcode`\@=11
\font\tenmsa=msam10
\font\sevenmsa=msam7
\font\fivemsa=msam5
\font\tenmsb=msbm10
\font\sevenmsb=msbm7
\font\fivemsb=msbm5
\newfam\msafam
\newfam\msbfam
\textfont\msafam=\tenmsa  \scriptfont\msafam=\sevenmsa
  \scriptscriptfont\msafam=\fivemsa
\textfont\msbfam=\tenmsb  \scriptfont\msbfam=\sevenmsb
  \scriptscriptfont\msbfam=\fivemsb
\def\hexnumber@#1{\ifcase#1 0\or1\or2\or3\or4\or5\or6\or7\or8\or9\or
      A\or B\or C\or D\or E\or F\fi }

\font\teneuf=eufm10
\font\seveneuf=eufm7
\font\fiveeuf=eufm5
\newfam\euffam
\textfont\euffam=\teneuf
\scriptfont\euffam=\seveneuf
\scriptscriptfont\euffam=\fiveeuf
\def\frak{\ifmmode\let\next\frak@\else
 \def\next{\errmessage{Use \string\frak\space only in math mode}}\fi\next}
\def\goth{\ifmmode\let\next\frak@\else
 \def\next{\errmessage{Use \string\goth\space only in math mode}}\fi\next}
\def\frak@#1{{\frak@@{#1}}}
\def\frak@@#1{\fam\euffam#1}

\edef\msa@{\hexnumber@\msafam}
\edef\msb@{\hexnumber@\msbfam}
\mathchardef\square="0\msa@03
\mathchardef\subsetneq="3\msb@28
\mathchardef\ltimes="2\msb@6E
\mathchardef\rtimes="2\msb@6F
\def\Bbb{\ifmmode\let\next\Bbb@\else
\def\next{\errmessage{Use \string\Bbb\space only in math mode}}\fi\next}
\def\Bbb@#1{{\Bbb@@{#1}}}
\def\Bbb@@#1{\fam\msbfam#1}
\catcode`\@=12
%
%


\def\PP{{\Bbb P}}
\def\F{{\bf F}}

\def\tto{\longrightarrow}

\def\today{\ifcase\month\or
 Jan\or Febr\or  Mar\or  Apr\or May\or Jun\or  Jul\or
 Aug\or  Sep\or  Oct\or Nov\or  Dec\or\fi
 \space\number\day, \number\year}
\vskip 6.5pc
\noindent
\font\eighteenbf=cmbx10 scaled\magstep2
\vskip 2.0pc
\centerline{\eighteenbf An Asymptotically Good Tower of Curves }
\smallskip
\centerline{\eighteenbf Over the Field  with Eight Elements}
\smallskip
\medskip
\vskip 2pc
\font\titlefont=cmcsc10 at 11pt
\centerline{\titlefont Gerard van der Geer and Marcel van der Vlugt}
\vskip 2.0pc
\bigskip
\centerline{\bf Introduction}
\medskip
In this note we construct an explicit asymptotically good tower
of curves over the field $\F_8$. For a curve $C$ defined over a
finite field $\F_q$ of cardinality
$q$ we denote by $\#C(\F_q)$ the number of $\F_q$-rational points on 
$C$. If furthermore we denote, as usual, 
by $N_q(g)$ the function
$$
N_q(g)= \max \{ \# C(\F_q): C/F_q, \, g(C)=g \},
$$
where $C$ runs through the set of smooth absolutely irreducible projective
curves of genus $g$
defined over $\F_q$ then Drinfeld and Vladuts showed in [D-V] the
inequality
$$
\lim\sup_{g \to \infty} {N_q(g) \over g} \leq \sqrt{q}-1, \eqno(1)
$$
and Ihara constructed in [I] for $q$ a square a sequence of modular curves
which attains the upper bound in (1).

It then came as a surprise when in  1995
Garcia and Stichtenoth constructed in [G-S1]
for $q$ a square a tower
of Artin-Schreier covers
$$
\ldots \tto  C_i \tto C_{i-1}\tto\ldots \tto C_1 \tto C_0
$$
which is defined over $\F_q$ and given by a simple recursive equation
such that
$$
\lim_{i\to \infty} g(C_i)=\infty \qquad {\rm and} \quad
  \lim_{i\to \infty} {\# C_i(\F_q) \over g(C_i)}=
\sqrt{q}-1.
$$

An infinite tower $C_{\bullet}$ of covers of curves over $\F_q$
$$
\ldots \tto C_i \tto C_{i-1}\tto\ldots \tto C_1 \tto C_0
$$
with $g(C_i)>1$ for some $i\geq 0$ is called an {\sl 
asymptotically good tower}
if 
$$
\ell(C_{\bullet})= \lim_{i \to \infty} { \# C_i(\F_q) \over g(C_i)} > 0.
$$
Note that in [G-S2] it is shown that this limit exists for 
towers having at least one index $i$ with $g(C_i)>1$.

Apart from having 
an evident charm of their own, asymptotically good towers
are important for coding theory, since such towers enable the construction 
of long error correcting codes over $\F_q$ which can correct a
fixed percentage of errors per codeword and have a positive 
transmission rate. However, for this application it is essential that the
curves are in explicit form.

For $q$ is not a square the results are much less complete. It is not
known how good the Drinfeld-Vladuts upper bound (1) is for that case.
For $q$ not a square asymptotically good towers of curves are mainly obtained
by class field theory, see for example [N-X]. These constructions are not 
explicit. In 1985 Zink using certain Shimura surfaces 
constructed in [Z] a tower $C_{\bullet}$ of curves 
defined over $\F_{p^3}$ for $p$ a prime with limit
$$
\ell(C_{\bullet})\geq { 2(p^2-1)\over p+2}, \eqno(2)
$$
but that construction is far from explicit.
In [G-S-T] there is an explicit asymptotically good tower of Kummer 
covers over $\F_{q=p^m}$ for $m\geq 2$ with limit
$$
\ell(C_{\bullet}) \geq { 2 \over q-2}.
$$

Here we present an explicit tower $C_{\bullet}$ 
of Artin-Schreier curves defined over $\F_8$
given by a simple recursive equation
with limit
$$
\ell(C_{\bullet})= 3/2.
$$
One should compare this with (2). It remains an interesting problem
to see whether our explicit tower is related to that of Zink, cf. the
remarks made by Elkies at the end of [E].
Another interesting problem is to extend our construction to other
fields of odd degree over the prime field.

We give explicit formulas for the genus and number of rational points for
the curves $C_i$ in our tower. The ramification behavior 
turns out to be rather subtle
with alternating ramification and non-ramification.

\bigskip

\centerline{\bf \S1 The basic equation}
\bigskip
\noindent
In our search for curves over finite fields with many points
we came across a curve defined over $\F_8$ with a remarkable
property. The curve of genus $1$ given by the affine equation
$$
x_1^2+x_1=x_0+1+1/x_0
$$
has $14$ $\F_8$-rational points and attains the Hasse-Weil bound
for $\F_8$. To each $x_0\in \F_8-\F_2$ there correspond $2$
solutions $x_1 \in \F_8-\F_2$ and if $x_0$ runs through $\F_8-\F_2$
then so does $x_1$. This implies that the system of equations
$$
x_{i+1}^2+x_{i+1}= x_i + 1 +1/x_i
\qquad i=0,1,2,\ldots
$$
has sequences of solutions $(x_0,x_1,x_2,\ldots)$ for
every $x_0 \in \F_8-\F_2$.
\smallskip
Consider in $\PP^1\times \PP^1$ over the field $\F_2$ 
the closure of the affine curve given by the
equation
$$
x_1^2+x_1= x_0+1+1/x_0.
$$
This defines a smooth projective 
curve $C$ of genus $1$ together with two morphisms
$ b_1: C \to \PP^1$, $(x_0,x_1) \mapsto x_0$,  and  
$e_1: C \to \PP^1$, $(x_0,x_1) \mapsto x_1$ of degree~$2$. 
The curve $C$ possesses $2$ points rational over $\F_2$,   
$8$  points rational over $\F_4$ 
and $14$ points rational over $\F_8$.  The correspondence $C$ 
in $\PP^1 \times \PP^1$  preserves the points of $\PP^1(\F_8)-\PP^1(\F_2)$
and surprisingly also those of $\PP^1(\F_4)$.
\bigskip
We consider the following infinite tower 
$C_{\bullet}$ of smooth projective 
curves defined over~$\F_2$
$$
\tto C_i {\buildrel \pi_i \over \tto} C_{i-1} {\buildrel \pi_{i-1} \over
\tto} \ldots {\buildrel\pi_2 \over
\tto} C_1 {\buildrel\pi_1 \over \tto} C_0= \PP^1,
$$
where we take an affine coordinate $x_0$ on $C_0$ and where the cover
$C_i \to C_{i-1}$ is given by the affine equation
$$
x_i^2+x_i= x_{i-1}+1+ {1 \over x_{i-1}}\qquad \hbox {\rm for $i\geq 1$}.
\eqno(3)
$$
Equivalently, we can describe the curve $C_i$ as the 
normalization of the curve $D_i$ defined by
$$
D_i= \{ (p_0,p_1,\ldots,p_{i})\in \PP^1 \times \ldots \times \PP^1
: (p_j,p_{j+1}) \in C \, \hbox {\rm for 
$j=0,\ldots,i-1$)} \}.
$$
This shows that $C_i$ for $i\geq 1$ admits two maps 
$b_i: C_i \tto C_{i-1} $ (resp.\ $e_i: C_i \tto C_{i-1}$)
given (on the model $D_i$) 
by  $(p_0,\ldots , p_i) \mapsto (p_0,\ldots,p_{i-1})$
(resp.\ $\mapsto (p_1,\ldots,p_i)$). The curve $C_i$ is then the normalization
of the fibre
product 
$$
C_{i-1} \times_{C_{i-2}} C_{i-1}
$$
via the maps $b_{i-1}$ and $e_{i-1}$.
\bigskip

We now work over the algebraic closure $\F$ of $\F_2$ and
consider geometric points of $C_i\otimes \F$.
It will turn out that ramification in $C_i/C_{i-1}$ can occur only
in points $P$ that map to a point of $D_i$ with coordinates
in $\PP^1(\F_4)$. We therefore introduce the following
notation for such points. 
By $P=P(a_0,a_1,\ldots,a_i)$ 
with $a_j \in \PP^1(\F_4)$ 
we denote a 
point on $C_i$ such that $x_j(P)=a_j$ for $0 \leq j \leq i$. That is,
the point $(a_0,\ldots, a_i)$ is the image point in $D_i$ and
will be called the index sequence of the point.
Note that  because of the normalization, a point $P(a_0,\ldots,a_i)$
of $C_i$ is not
necessarily uniquely determined by 
its  index sequence $(a_0,\ldots,a_i)$ on $D_i$.

We shall write $\rho$ for a primitive element of $\F_4$.
Note that in an index sequence $(a_0,a_1,\ldots,a_i)$ of a point $P$
we have
$$
\matrix{\infty\cr 0 \cr 1 \cr \rho \cr \rho^2\cr}\qquad \hbox{\rm is
followed by} \qquad\matrix{\infty \cr \infty \cr \hbox{\rm $\rho$ or $\rho^2$}
\cr \hbox{\rm $0$ or $1$}\cr \hbox{\rm $0$ or $1$} \cr}
$$
Sometimes we shall write $(a_0,\ldots,a_i,\infty^j)$ for
a point $(a_0,\ldots,a_i,\underbrace{\infty, \ldots, \infty}_{j \times})$.
\bigskip
\centerline{\bf \S 2 The principal part of the $x_i$}
\bigskip
\noindent
The main problem to find the limit $\ell(C_{\bullet})$
of our tower lies in  the determination of the genus $g(C_{i})$.
In order to compute it  we have to find the ramification
divisor of $C_{i+1}$ over $C_i$. 
Since we are dealing with  Artin-Schreier equations we can restrict
ourselves to the points which are poles of the function $f_i= x_i+1+1/x_i$.
The contribution to the ramification 
is determined by the orders ${\rm ord}_P(f_{i}^*)$   for the poles $P$ on
$C_{i}$ of the Artin-Schreier reduction $f_i^*$ of the 
function $f_{i}$.

\proclaim (2.1) Lemma. The zeros of $x_i$ on $C_i$ are of the form 
$P(a_0,a_1,\ldots,a_i)$
with $a_i=0$, $a_{i-j}\in {\F}_4-{\F}_2$ for $j\geq 1$  odd and $a_{i-j}=1$ for
$j\geq 2$ even. The poles of $x_i$ are of the form
 $P(b_0,b_1,\ldots,b_j,\infty^{i-j})$ with
$0\leq j\leq i-1$ and $(b_0,\ldots,b_j)$ an index sequence of
a zero of $x_j$  or of the form $P(\infty^{i+1})$.
\par
\noindent
{\sl Proof.} By induction on $i$. 
The lemma is true for $x_0$. From the equation
$$
x_i^2+x_i= {x_{i-1}^2 +x_{i-1}+1\over x_{i-1}} =f_{i-1}
$$
it follows that we have the equality of divisors on $C_i$
$$
(f_{i-1})= (x_i)+(x_i+1)=(x_i)_0 + (x_i)_1 -2(x_i)_{\infty}.
$$
So poles of $x_i$ lie above poles of $f_{i-1}$ and the points $P$ on $C_i$
with $x_i(P) \in \F_2$ lie above the zeros of $f_{i-1}$.
Moreover, we have
$$
\eqalign{
(f_{i-1})&= (x_{i-1}+\rho)+(x_{i-1}+\rho^2)-(x_{i-1})\cr
&= (x_{i-1})_{\rho} + (x_{i-1})_{\rho^2} -(x_{i-1})_0-(x_{i-1})_{\infty}.\cr
}
$$
which implies that the poles of $f_{i-1}$ are the zeros and poles of $x_{i-1}$,
while the zeros of $f_{i-1}$ are the points $P$ on $C_{i-1}$ with
$x_{i-1}(P)\in {\F}_4 - {\F}_2$. Hence the poles of $x_i$ are the points on $C_i
$
above the zeros and the poles of $x_{i-1}$ on $C_{i-1}$, 
whereas the zeros of $x_i$ lie above
 points $P$ on $C_{i-1}$
with $x_{i-1}(P) \in {\F}_4 - {\F}_2$. So we obtain the index sequence of a pole
of $x_i$ by adding $\infty$ to a  zero  or pole of $x_{i-1}$
and we obtain the index sequence of a zero of $x_i$ by adding a zero
to an index sequence which ends with an element of ${\F}_4 - {\F}_2$
and in which $1$ and elements of $\{ \rho, \rho^2\}$ alternate. $\square$
\bigskip
In the following we shall develop rational functions on $C_i$
as a power series 
in a local parameter at a given point $P$, that is, we consider
the function as an element of the quotient field of the completion
of the local ring of $P$. Often we are only interested in the
principal part and neglect elements that are regular,
i.e. elements of  (the completion of) the local ring.
By the notation
$$
f=g+O(P)
$$
we mean that $f-g$ is regular in $P$, that is, is an element of
$O_P$ or of $\hat{O}_P$.
\bigskip
Consider now a sequence  of points $P_0 \in C_0, P_1 \in C_1, \ldots,
P_i \in C_i$ with $\pi_{\ell}(P_{\ell})=P_{\ell-1}$ for $\ell=1,\ldots,i$
and with the property that $1$ and $\rho$ or $\rho^2$ alternate in the
index sequence $(a_0,\ldots,a_i)$ of $P_i$.

We shall first assume that $a_0=1$. Then the function $t=x_0+1$
provides a local parameter at $P_0$ on $C_0$ 
and the pull back (under the maps $\pi_{\ell}$) of this function
(again denoted by $t$) is still a local parameter at the points $P_j$
on $C_j$ for $j\leq i$.

In the completion of the local ring $\hat{O}_{P_j}\cong \F[[t]]$ the function $x_j$
can be written as a power series in $t$
$$
x_j= a_j +m_j(t),
$$
where $m_j(t) \in \F[[t]]$ has ${\rm ord}_t(m_j)\geq 1$.

\proclaim (2.2) Lemma. In the quotient field $\F((t))$ of the
formal power series ring $\hat{O}_{P_j}\cong \F[[t]]$ 
 the function $m_j(t)$ satisfies for $0\leq j \leq i$
the relations
$$
\eqalign{
{1 \over m_j}& = {a_{j-1} \over m_{j-1}}+O(P_j) 
 \qquad \hbox{\rm for $j \geq 2$ even},\cr
{1 \over m_j}& = {1 \over m_{j-1}^2}+ { 1 \over m_{j-1}}  +O(P_j)
 \qquad \hbox{\rm for $j $ odd.}\cr
}
$$
\par
\smallskip
\noindent
{\sl Proof.} We start with $m_0(t)=t$. For even $j\geq 2$ we have 
$a_{j-1} \in \{ \rho, \rho^2\}$ and $a_j=1$ since we assumed that 
$a_0=1$. From the relation
$$
x_j^2+x_j = x_{j-1}+ 1 + 1/x_{j-1}
$$
we obtain
$$
\eqalign{
m_j^2+m_j &= a_{j-1}+m_{j-1} +1 + 1/(a_{j-1}+m_{j-1})\cr
&= a_{j-1} + m_{j-1} + 1 + (1 /a_{j-1})
\sum_{n=0}^\infty (m_{j-1} / a_{j-1})^n\cr
&= a_{j-1}^2m_{j-1} + m_{j-1}^2 +
\hbox{\rm higher powers of $m_{j-1}$}.\cr
}
$$
This implies that $m_j$ is the product of $a_{j-1}^2m_{j-1}$ with
a $1$-unit $u$ in $m_{j-1}$, i.e.\ a unit 
of the form $u=1+r$ with $r\in (m_{j-1})$. 
So we get
$$
{1 \over m_{j}}= {a_{j-1} \over m_{j-1}} \cdot u =
{a_{j-1} \over m_{j-1}} + O(P_{j}).
$$
For $j$ odd we have $a_j \in \{ \rho, \rho^2\}$ and $a_{j-1}=1$. In the
same way as for $j$ even we obtain
$$
a_j^2+m_j^2 +a_j+m_j= 1 + m_{j-1} +1 + {1 \over 1 + m_{j-1}}=
1 + \sum_{n=2}^{\infty} m_{j-1}^n,
$$
that is,
$
m_j^2+m_j= \sum_{n=2}^{\infty} m_{j-1}^n
$,
so that  we have
$$
m_j= m_{j-1}^2+m_{j-1}^3 +\hbox{\rm higher powers of $m_{j-1}$}.
$$
This means that
$$
\eqalign{
{1 \over m_{j}} &= {1 \over m_{j-1}^2} + { 1 \over m_{j-1}} +
\hbox{\rm higher powers of $m_{j-1}$}\cr
&= {1 \over m_{j-1}^2} + {1 \over m_{j-1}} +
O(P_j).
}
$$
This completes the proof of the lemma. $\square$
\bigskip

We denote the principal part of $1/m_j$ by $F_j$. We now can deduce the
following corollary.

\smallskip
\noindent
\proclaim (2.3)  Corollary. The principal part $F_j$ of $1/m_j$ satisfies:
$$
F_j=\cases{
F_{j-1}^2+F_{j-1} & for $j$ odd\cr
a_{j-1}\cdot F_{j-1} & for $j\geq 2$ even.\cr}
$$
Furthermore, $F_j$ is a $2$-linearized polynomial in $1/t$ of the form
$$
F_j= {b_{k} \over t^{2^k}} + {b_{k-1} \over t^{2^{k-1}}} + \ldots + 
{b_0 \over t},
$$
where $k=[(j+1)/2]$,   the coefficients $b_{\ell}$ are in $\F_4$ and
$b_{k}\neq 0$. \par

\smallskip
\noindent
{\sl Proof.} The relations for $F_j$ follow at once from Lemma.(2.2). We
have $F_0=1/t$ and $F_1= 1/t^2+1/t$ from which the
formula for $F_j$ follows by induction. $\square$
\bigskip

For an index sequence $(a_0,a_1,\ldots, a_i)$ where $a_0 \in \{ \rho, \rho^2\}$
we have a similar result.
\bigskip
\centerline{\bf \S3 The ramification behavior}
\bigskip
Now we study the ramification behavior in a  point
$P_i= P(a_0,\ldots,a_{i-1},a_i=0)$ which is a zero of $x_i$ on $C_i$ for
$i\geq 2$. We assume $a_0=1$. Then $a_{\rm odd}\in \{ \rho, \rho^2\}$
and $i$ is even. In a point $P_i$ where $a_0 \in \{ \rho, \rho^2\}
$ the ramification behavior is similar. 

Since we are working with Artin-Schreier covers here we introduce the
standard notation $\wp(f)= f^2+f$ for an element $f$ 
in one of our function fields.

\proclaim (3.1) Lemma. A linear combination $\sum_{j=2, \, \rm even}^i
B_{j,i} F_j$  with coefficients $B_{j,i}\in \F_4$ can be written as
$$
\sum_{j=2, \, \rm even}^i B_{j,i} F_j= \wp(\sum_{j=0, \, \rm even}^{i-2}
B_{j,i-2}F_j) + B_i^* F_0 \eqno(4)
$$
with
$$
B_i^*= \wp(\sum_{j=2,\, \rm even}^i B_{j,i}a_{j-1}  
) \eqno(5)
$$
and
$$
B_{j,i-2}=\big( B_i^*+\wp( \sum_{k=2,\, \rm even}^j B_{k,i}a_{k-1}) 
\big)a_{j+1}^2 +
\wp(B_{j+2,i}). \eqno(6)
$$
\par
\smallskip
\noindent
{\sl Proof.} Using Corollary (2.3)  we find for even $j \geq 2$
$$
B_{j,i}F_j = B_{j,i}a_{j-1}F_{j-1} = B_{j,i}a_{j-1}\wp(
F_{j-2})=
\wp(B_{j,i}^2a_{j-1}^2 F_{j-2})+\wp(B_{j,i}a_{j-1})F_{j-2},
$$
that is
$$
B_{j,i}F_j= 
\wp(B_{j,i}^2a_{j-1}^2 F_{j-2})+\wp(B_{j,i}a_{j-1})F_{j-2}.
\eqno(7)
$$
Applying (7) to the second term in the RHS of (7) we obtain
$$
\wp(B_{j,i}a_{j-1})F_{j-2} = \wp((
\wp(B_{j,i}a_{j-1})a_{j-3}^2F_{j-4})+
\wp(B_{j,i} a_{j-1})F_{j-4},
$$
where we use $a_{\rm odd}^2+a_{\rm odd}=1$.
Continuing this way we find an expression
for $B_{j,i}F_j$ as 
$$
\wp( \hbox{\rm linear combination of $F_{j-2},\ldots,
F_0$ })
+\wp(B_{j,i}a_{j-1})F_0.
$$
Adding these relations for all terms in $\sum_{j=2, \, \rm even}^i B_{j,i}F_j$ we find
(4) with coefficients satisfying the equations (5) and (6). $\square$

\bigskip
Note that all coefficients are in $\F_4$ and that $B_i^*$ is in $\F_2$.

\bigskip
The cover $C_{i+1}/C_{i}$ is given by the equation
$$
x_{i+1}^2+x_{i+1}= x_i + 1 +1/x_i.
$$
In a point $P_i $ with index sequence $(1,a_1,\ldots,a_{i-1},0)$
we have the relation
$$
x_{i+1}^2+x_{i+1}  = {1 \over m_i} + O(P_i) = F_i+O(P_i).    \eqno(8)
$$
Therefore the principal part of 
$x_{i+1}^2+x_{i+1}$ in $P_i$ is $F_i$.
According to Lemma  (3.1) we have
$$
F_i= \wp(\sum_{j=0, \, \rm even}^{i-2} B_{j,i-2} F_j) + B_i^*F_0
$$
with $B_i^*= a_{i-1}^2+a_{i-1}=1$ and  $B_{j,i-2}=a_{j+1}^2$ for even
$j=0, 2, \ldots, i-2$.

By the substitution
$$
X_{i+1}=  \sum_{j=0}^{i-2} B_{j,i-2} F_j + x_{i+1}
$$
we can reduce the equation (8) to 
$$
X_{i+1}^2+X_{i+1}= B_i^*F_0 + O(P_i) 
= F_0 + O(P_i)
 \eqno(9)
$$
with $F_0=1/t$.

\proclaim (3.2) Corollary. The point $P_i=P(a_0=1,a_1,\ldots,a_i=0)$
is totally ramified in $C_{i+1}/C_i$ and the contribution of $P_i$
to the ramification divisor of $C_{i+1}/C_i$ is $2$. In the pole $P_{i+1}
=P_i(\infty)$
of $x_{i+1}$ we have
$$
{\rm ord}_{P_{i+1}}(x_{i+1})= -2^{[(i+1)/2]}.
$$
\par
\bigskip
As the next step we consider the behavior of the pole $P_{i+1}=P_i(\infty)$
of the function 
$f_{i+1}=x_{i+1} + 1 + (1/x_{i+1})$ in the cover $C_{i+2}/C_{i+1}$.

In $P_{i+1}$ the equation of $C_{i+2}/C_{i+1}$ is
$$
\eqalign{
x_{i+2}^2+x_{i+2}&= x_{i+1} + O(P_{i+1}) \cr
&=\sum_{j=0,\, \rm even}^{i-2} B_{j,i-2}F_j + X_{i+1} + O(P_{i+1}).\cr
} \eqno(10)
$$
If we apply Lemma (3.1)
  to the linear combination $\sum_{j=2,\, \rm even}^{i-2} B_{j,i-2} F_j$
the RHS of (10) becomes
$$
\wp(\sum_{j=0, \, \rm even}^{i-4} B_{j,i-4}F_j) +
B_{i-2}^* F_0 + B_{0,i-2}F_0 + X_{i+1} + O(P_{i+1}).
\eqno(11)
$$
Then by using (9),  i.e., by
substituting $F_0= X_{i+1}^2+X_{i+1} + O(P_{i+1})$,
the expression (11) is converted to
$$
\wp(\sum_{j=0, \, \rm even}^{i-4} B_{j,i-4} F_j)+ 
\wp(B_{i-2}^* X_{i+1}) + \wp (B_{0,i-2}^2X_{i+1})+
(B_{0,i-2}^2+B_{0,i-2}+1)X_{i+1} + O(P_{i+1}) 
$$
with $B_{0,i-2}^2+B_{0,i-2}+1= a_1^2+a_1+1=0$.
Hence the equation of $C_{i+2}/C_{i+1}$ in $P_{i+1}$ is of the
form
$$
x_{i+2}^2+x_{i+2} = \wp( \sum_{j=0}^{i-4} B_{j,i-4}F_j) +
\wp((B_{i-2}^* + B_{0,i-2}^2)X_{i+1}) +
O(P_{i+1}).\eqno(12)
$$

\proclaim (3.3) Corollary. The pole $P_{i+1}$ of $f_{i+1}$  is unramified in the cover
$C_{i+2}/C_{i+1}$ and in a point $P_{i+2}=P_{i+1}(\infty)$ above $P_{i+1}$
we have
$$
{\rm ord}_{P_{i+2}}(x_{i+2})= -2^{[(i+1)/2]-1}.
$$
\par

Note that (12) implies that
$$
x_{i+2}+(\sum_{j=0, \, \rm even}^{i-4} B_{j,i-4}F_j) +
(B_{i-2}^*+B_{0,i-2}^2)X_{i+1}
$$
is integral in the point $P_{i+2}$. 

To analyze the situation in $P_{i+2}$ we start with the equation of $C_{i+3}/
C_{i+2}$ in this point:
$$
\eqalign{
x_{i+3}^2+x_{i+3}&= x_{i+2} + 
O(P_{i+2})\cr
&=\big(\sum_{j=0,\, \rm even}^{i-4} B_{j,i-4}F_j\big)+
(B_{i-2}^* +B_{0,i-2}^2)X_{i+1} +
O(P_{i+2}).  \cr}\eqno(13)
$$
Using Lemma  (3.1) and (9) the RHS of (13) is of the form
$$
\eqalign{
\wp(\sum_{j=0, \, \rm even}^{i-6} B_{j,i-6}F_j) + &
\wp((B_{i-4}^*+B_{0,i-4}^2) X_{i+1})+ \cr
+ (\wp(B_{0,i-4})+ & B_{i-2}^*+B_{0,i-2}^2) X_{i+1} + O(P_{i+2}).
\cr} \eqno(14)
$$
Since (6) implies $\wp(B_{0,i-4})=B_{i-2}^*$ the coefficient 
of $X_{i+1}$ in (14) is $B_{0,i-2}^2=a_1^2$. So the right  hand side
of (13) has the form
$$
\wp(\gamma) + B_{0,i-2}^2X_{i+1} + O(P_{i+2})
$$
for some $\gamma \in \F(C_{i+2})$.
If we set 
$$
X_{i+3}= \sum_{j=0}^{i-6} B_{j,i-6}F_j +
(B_{i-4}^*+B_{0,i-4}^2)X_{i+1}
$$
the equation of $C_{i+3}/C_{i+2}$ becomes
$$
X_{i+3}^2+X_{i+3}= B_{0,i-2}^2 X_{i+1} +
O(P_{i+2}).
\eqno(15)
$$

\proclaim (3.4) Corollary. 
The point $P_{i+2}$ is totally ramified in the cover $C_{i+3}/C_{i+2}$
and the contribution to the ramification divisor is $2$.
In $P_{i+3}=P_{i+2}(\infty)$ above $P_{i+2}$ we have
$$
{\rm ord}_{P_{i+3}}(x_{i+3})= -2^{[(i+1)/2]-1}.
$$
\par

If we continue along these lines we obtain the following formulas:

\proclaim (3.5) Formula. 
For $t$ even and $2\leq t \leq i$ the equation of $C_{i+t}/C_{i+t-1}$
in a point $P_{i+t-1}$ is
$$
\eqalign{
x_{i+t}^2+x_{i+t}=\wp(\sum_{j=0,\, \rm even}^{i-2t} B_{j,i-2t}F_j) +
\wp((B_{i-2t+2}^*+B_{0,i-2t+2}^2)X_{i+1})+ &\cr
 \quad +\sum_{k=1}^{t/2-1} \wp((B_{2k-2,i-2t+4k}B_{2k-2,i-2}^2)X_{i+2k+1})
+& O(P_{i+t-1}). \cr
}
$$
\par

\proclaim (3.6) Formula. 
For $t$ odd and $3\leq t \leq i-1$ the equation of $C_{i+t}/C_{i+t-1}$
in a point $P_{i+t-1}$ is
$$
\eqalign{
x_{i+t}^2+x_{i+t}=\wp(\sum_{j=0, \, \rm even}^{i-2t} B_{j,i-2t}F_j) +
\wp((B_{i-2t+2}^*+B_{0,i-2t+2}^2)X_{i+1})+ &\cr
 \quad+\sum_{k=1}^{(t-3)/2} \wp((B_{2k-2,i-2t+4k}B_{2k-2,i-2}^2)X_{i+2k+1})
+B_{t-3,i-2}^2X_{i+t-2} +&O(P_{i+t-1}). \cr
}
$$
\par
\smallskip
\noindent
{\bf (3.7) Remark.} 
The function $X_{i+2k+1}$ with $k\geq 1$ satisfies an equation of the form
$$
X_{i+2k+1}^2+X_{i+2k+1}=B_{2k-2,i-2}^2 X_{i+2k-1} +
O(P_{i+2k}).
$$
We also find
$$
{\rm ord}_{P_{i+t}}(x_{i+t})=-2^{[(i+2)/2]-[t/2]} \qquad {\rm for} \quad
2 \leq t \leq i.
$$
Hence ${\rm ord}_{P_{2i}}(x_{2i})=-1$ and this implies that from $P_{2i}$
on the extensions in the tower are totally ramified above $P_{2i}$
and their contribution to the ramification divisor is $2$.
\bigskip
We summarize the preceding results in the following theorem.

\proclaim (3.8) Theorem. A pole $P_{i+j}=P(a_0,\ldots,a_{i-1},0,\infty^j)$
of $x_{i+j}$ on $C_{i+j}$ for $j\geq 1$ (resp.\ of $1/x_i$ for $j=0$)
with $a_0=1$ is
\item{1)} totally ramified in $C_{i+j+1}/C_{i+j}$ for
$j$ with $j=0,2,4,\ldots,i-2$ or with $j\geq i$
and each of these contributes  $2$ to the ramification divisor,
\item{2)} unramified in $C_{i+j+1}/C_{i+j}$ for $j=1,3,5,\ldots,i-1$.\par

For a pole whose index sequence starts with $a_0 \in \{ \rho, \rho^2\}$
there is the following  similar result.

\proclaim (3.9) Theorem. A pole $P_{i+j}=P(a_0,\ldots,a_{i-1},0,\infty^j)$
of the function $x_{i+j}$ on $C_{i+j}$
with $a_0\in \{\rho, \rho^2\}$ is
\item{1)} totally ramified in $C_{i+j+1}/C_{i+j}$ for all $j$ with
$j=0,2,4,\ldots,i-3$
 and for all $j$ with $j\geq i-1$
with contribution $2$ to the ramification divisor,
\item{2)} unramified in $C_{i+j+1}/C_{i+j}$ for $j=1,3,5,\ldots,i-2$.\par
\noindent
{\bf (3.10) Remark.} 
We always have the totally ramified points $P(\infty,\infty,
\ldots,\infty)$
and also $P(0,\infty,\infty,\ldots,\infty)$ 
which contribute $2$ to the ramification divisor.

\bigskip
\centerline{\bf \S4 The  genus and the number of points in the tower}
\bigskip
In order to compute the genus of our curve $C_i$ 
we have to count the number of points on the curve $C_i$ which
contribute to the ramification divisor.
We find using (3.8)--(3.10)

\proclaim (4.1) Theorem. Let $n_i$ be the number of points on $C_i$ which
are totally ramified in $C_{i+1}/C_i$. Then
$$
n_i=\cases{ ([(i+2)/ 4]+2) 2^{i/2} & for $i$ even,\cr
([i/ 4]+2) 2^{(i+1)/2} & for $i$ odd.\cr
}
$$
\par

Now it is not difficult to determine the genus $g(C_i)$ of $C_i$.
From the Hurwitz formula it follows that
$$
g(C_i)=1+ \sum_{j=1}^{i-1} 2^{i-j-1} n_j. \eqno(16)
$$
If we combine (16) with  theorem (4.1) we get:

\proclaim (4.2) Theorem. The genus $g(C_i)$ of $C_i$ satisfies
$$
g(C_i)= 2^{i+2}+1-\cases{ (i+10)2^{(i/2)-1}  & for $i$ even, \cr
(i+2[i / 4]+15) 2^{(i-3)/2}  & for $i$ odd.\cr
}
$$

Now we count the number $\#C_i(\F_8)$
 of $\F_8$-rational points on $C_i$. 

\proclaim (4.3) Theorem. We have $\# C_i(\F_8)= 6\cdot 2^i +2$.
\par

\noindent
{\sl Proof.} Let $\alpha$ be a primitive element of $\F_8$ which
satisfies $\alpha^3+\alpha +1=0$. For $x \in \F_8-\F_2$ we find
$$
\{ x+{1 \over x} +1 : x \in \F_8-\F_2\} = \{ \alpha, \alpha^2,\alpha^4\},
$$
but also
$$
\{ y^2+y: y \in \F_8-\F_2\} =\{ \alpha, \alpha^2,\alpha^4\}.
$$
This means that a point $x \in \PP^1(\F_8)$ with $x \not\in \PP^1(\F_2)$
splits completely in the tower. This yields $6 \cdot 2^i$ rational
points over $\F_8$ on $C_i$. Besides these we have two totally
ramified points $P(0,\infty,\ldots,\infty)$ and $P(\infty, \ldots, \infty)$
defined over $\F_2$.
$\square$
\smallskip

Combining this with the formula for the genus we obtain the following Theorem.

\proclaim (4.4) Theorem. The tower of curves $C_{\bullet}$ over $\F_8$ is asymptotically good with limit
$$
\lim_{i\to \infty} { \# C_i(\F_8) \over g(C_i)} = {3 \over 2}.
$$
\par
\bigskip
\centerline{\bf References}
\bigskip
\noindent
[D-V] V.G.\ Drinfeld, S.G. Vladuts: Number of points of an algebraic curve.
{\sl Funct.\ Anal.\ \bf 17} (1983), p.\ 68--69.
\smallskip
\noindent
[E] N.\ Elkies: Explicit modular towers. Preprint Harvard University, 2000.

\smallskip
\noindent
[G-S1] A.\ Garcia, H.\ Stichtenoth: A tower of Artin-Schreier extensions of
function fields attaining the Drinfeld-Vladut bound.
{\sl Inv.\ Math.\ \bf 121} (1995), p.\ 211--222.

\smallskip
\noindent
[G-S2] A.\ Garcia, H.\ Stichtenoth: 
On the asymptotic behavior of some towers of function fields over
finite fields. {\sl J.\ of Number Th. \bf 61} (1996), p.\ 248--273.

\smallskip
\noindent
[G-S-T] A.\ Garcia, H.\ Stichtenoth, M.\ Thomas:
 On towers and composita of towers of function fields over finite fields.
{\sl Finite Fields and their Appl.\ \bf 3} (1997), p.\ 257--274.

\smallskip
\noindent
[I] Y.\ Ihara: Some remarks on the number of points of algebraic curves over finite fields. {\sl J.\ Fac.\ Sci.\ Tokyo, Ser.\ Ia, \bf 28} (1982), p.\ 721--724.

\smallskip
\noindent
[N-X] H.\ Niederreiter, C.\ Xing: Global function fields with many
rational places and their applications.
{\sl Contemp.\ Math.\ \bf 225} (1999), p.\ 87--111.

\smallskip
\noindent
[Z] Th.\ Zink: Degeneration of Shimura curves and a problem in coding
theory. In {\sl Fundamentals of Computation Theory}, Springer LNCS 199,
p.\ 503--511. Springer, Berlin 1985.
\bigskip
\bigskip
\settabs3 \columns
\+G. van der Geer  &&M. van der Vlugt\cr
\+Faculteit
Wiskunde en Informatica &&Mathematisch Instituut\cr
\+Universiteit van
Amsterdam &&Rijksuniversiteit te Leiden \cr
\+Plantage Muidergracht 24&&Niels Bohrweg 1 \cr
\+1018 TV Amsterdam
&&2333 CA Leiden \cr
\+The Netherlands &&The Netherlands \cr
\+{\tt geer@science.uva.nl} &&{\tt vlugt@math.leidenuniv.nl} \cr

\bye